\documentclass[11pt, twosided]{article}
\usepackage{amsmath}
\usepackage{amsfonts}
\usepackage{color}
\usepackage{amsmath,amsthm,amssymb}
\makeatletter
\newtheorem{thm}{Theorem}[section]
\newtheorem{prop}[thm]{Proposition}
\newtheorem{lem}[thm]{Lemma}
\newtheorem{coro}[thm]{Corollary}
\theoremstyle{definition}
\newtheorem{defi}[thm]{Definition}

\newtheorem{exam}[thm]{Example}

\newtheorem{remark}[thm]{Remark}

\def\proof {{\noindent\bf Proof.}\quad}
\def \endproof{\hfill$\Box$\vspace {3mm}}

\newcommand{\ba}{\begin{array}}
\newcommand{\ea}{\end{array}}
\newcommand{\bt}{\begin{tabular}}
\newcommand{\et}{\end{tabular}}
\newcommand{\btb}{\begin{table}}
\newcommand{\etb}{\end{table}}
\newcommand{\bc}{\begin{center}}
\newcommand{\ec}{\end{center}}
\newcommand{\bea}{\begin{eqnarray}}
\newcommand{\eea}{\end{eqnarray}}
\newcommand{\Bea}{\begin{eqnarray*}}
\newcommand{\Eea}{\end{eqnarray*}}
\newcommand{\beq}{\begin{equation}}
\newcommand{\eeq}{\end{equation}}
\newcommand{\Beq}{\begin{equation*}}
\newcommand{\Eeq}{\end{equation*}}

\newcounter{vcenterstest}
\newlength{\leftlength}
\newlength{\rightlength}

\newlength{\vcentersskip}

\newcommand{\leftcentersright}[4][2]{%

\settowidth{\leftlength}{#2}%
        \settowidth{\rightlength}{#4}%

 \setcounter{vcenterstest}{#1}%
        \ifthenelse{\value{vcenterstest} = 0}

 {\setlength{\vcentersskip}{0pt}}{}%
        \ifthenelse{\value{vcenterstest} = 1}

 {\setlength{\vcentersskip}{\smallskipamount}}{}%

    \ifthenelse{\value{vcenterstest} = 2}

{\setlength{\vcentersskip}{\medskipamount}}{}%

   \ifthenelse{\value{vcenterstest} = 3}

{\setlength{\vcentersskip}{\bigskipamount}}{}%

 \ifthenelse{\value{vcenterstest} = 4}

{\setlength{\vcentersskip}{1cm}}{}%

\vskip\vcentersskip
        \noindent#2\hskip-\leftlength%

\hfill#3\hfill
        \mbox{}\hskip-\rightlength#4%

\vskip\vcentersskip
        }

\begin{document}
\title{Quasi-derivations  of Lie-Yamaguti algebras}
\author{Jie Lin$^{a}$, Yao Ma$^{b}$,
 Liangyun Chen$^{b}$ \thanks{Corresponding author.
{\em E-mail address:} chenly640@nenu.edu.cn.}\\
{\footnotesize\em $^{a}$ Civil Aviation University of China,
 Tianjin 300300, China}\\
{\footnotesize\em $^b$ Northeast Normal University, Changchun 130024,
China}}
\date{}

\maketitle

\baselineskip 18pt
\begin{abstract}
The concept of derivation for Lie-Yamaguti algebras is generalized in this paper. A quasi-derivation of an LY-algebra is embedded as derivation in a larger LY-algebra.  The relationship between quasi-derivations and robustness of Lie-Yamaguti algebras has been studied.
\medskip

\end{abstract}

\renewcommand{\theequation}{\arabic{section}.\arabic{equation}}

\section{Introduction}
A Lie-Yamaguti algebra is a binary-ternary algebra system which is denoted by $(T, \mu_1, \mu_2)$ (or $(T, [\cdot, \cdot], \{\cdot, \cdot, \cdot\})$) in this paper, and briefly called an LY-algebra(concretely, see Definition \ref{def1}). The LY-algebras with the binary multiplication $\mu_1=0$ are exactly the Lie triple systems, closely related with symmetric spaces, while the LY-algebras with the ternary multiplication $\mu_2=0$ are the Lie algebras. Therefore, they can be considered   as a simultaneous generalization of Lie triple systems and Lie algebras. They have been called ``generalized Lie triple systems¡° by Yamaguti in [\ref{ref01}](1957/1958) in an algebraic study of the characteristic properties of the
torsion and curvature of a homogeneous space with canonical
connection (the Nomizu¡¯s connection) [2](1954). Later on, these non-associative binary-ternary structures were called
¡°Lie triple algebras¡± by Kikkawa in [3](1975), then he studied the Killing-Ricci forms and invariant forms of Lie triple algebras respectively in [4](1981) and [5](1982).  The terminology of ¡°Lie-Yamaguti
algebras¡± is introduced by Kinyon and Weinstein in [6](2001) for these algebras.

LY-algebras have been treated by several authors in connection with geometric problems on homogeneous spaces ([\ref{ref04}]-[\ref{ref05}], [\ref{ref061}]-[\ref{ref063}]). Their structure theory has been studied by P. Bentio, C. Draper and A. Elduque in [\ref{ref07}]-[\ref{ref09}].% and their deformation theory is studied respectively by J. Lin, L. Y. Chen, Y. Ma, and T. Zhang, J. Li in [\ref{ref16}] and [\ref{ref17}].
 %Less known examples can be found in [\ref{ref010}] where a detailed analysis on the algebraic structure of LYAs arising from homogeneous spaces which are quotients of the compact Lie group $G_2$ is given.

Leger and Luks [\ref{ref18}] introduced a new concept called quasi-derivations of Lie algebras.
Inspired by them, we are interested in generalizing the derivations of LY-algebras. In this paper, we give the definitions and basic properties for generalized derivations, quasi-derivations, centroids, quasi-centroids of LY-algebars in section 2 and 3. In section 4, a quasi-derivation of an LY-algebra is embedded  as a derivation of a larger LY-algebra. Section 5 is devoted to the study of the connection of  quasi-derivations and robustness for LY-algebras, where the cohomology theory developed by Yamaguti ([\ref{ref091}]) is an important tool.

\section{Definitions and Notations }
\setcounter{equation}{0}
In this paper, $\mathbb K$ denotes a field.
\begin{defi}$^{[\ref{ref06}]}$ \label{def1}  A Lie-Yamaguti algebra(LY-algebra for short) is a vector space $T$
over $\mathbb K$ with a bilinear composition $[\cdot, \cdot]$ and a trilinear composition $\{ \cdot,  \cdot,
 \cdot\}$ satisfying: $$[a, a]=0, \eqno(LY1)$$
$$\{a, a, b\}=0, \eqno(LY2)$$ $$\{a, b, c\}+ \{b,c,a\}+\{c,a,b\}+[[a,b],c]+[[b,c],a]+[[c,a],b]=0, \eqno(LY3)$$ $$ \{[a,b], c, d\}+\{[b,c], a, d\}+\{[c,a], b, d\}=0,\eqno(LY4)$$
$$ \{a, b, [c,d]\}=[\{a, b, c\},d]+[c,\{a, b, d\}], \eqno(LY5)$$
$$\{a,b,\{c,d,e\}\}=\{\{a,b,c\},d,e\}+\{c,\{a,b,d\},e\}+\{c,d,\{a,b,e\}\}, \eqno(LY6)$$
for any $a, b, c, d, e \in T$.
\end{defi}
On a left Leibniz algebra $(L, \cdot)$, if we define $[x, y]:=\frac 12(x\cdot y-y\cdot x) $ (skew-symmetrization) and $\{x, y, z\}:=-\frac 14(x\cdot y)\cdot z$, then $(L, [\cdot, \cdot], \{\cdot, \cdot, \cdot\})$ is an LY-algebra([\ref{ref06}]).

${\rm End}(T)$ denotes the set of all linear maps of $T$.  Obviously, ${\rm End}(T)$ is a Lie algebra over $\mathbb K$ respect to the operation: $(D_1, D_2)\mapsto [D_1, D_2]=D_1D_2-D_2D_1.$ In order to distinguish the Lie structrue from the associative one, we write $\mathfrak {gl}(T)$ for ${\rm End}(T)$.

\begin{defi}[Ideal]
A subspace $I$ of an LY-algebra $T$ is called an ideal $T$ if $$ (x\in T, y\in I) \Longrightarrow [x, y]\in I \, \mbox{and}\, (x\in T, y\in T, z\in T)\Longrightarrow \{x, y, z\}\in T.$$
\end{defi}

\begin{defi}[Center]
The set $Z_T(I)=\{x\in T|\{x, a, y\}=\{y, a, x\}=0 \, \mbox{and}\, [x, y]=0, \forall a\in I, \forall y\in T\}$  is called the centralizer of $I$ in $T$. Particularly, $Z_T(T)$ denoted simply by $Z(T)$. We say that $T$ is centerless if $Z(T)=\{0\}.$
\end{defi}
\begin{defi}
The subalgebra $[T, T]+\{T, T, T\}$ of $T$ is called the derived algebra of $T$, we denote it by $T^{(1)}$.
\end{defi}
\begin{defi}[Derivation]
Let $(T, [\cdot, \cdot], \{\cdot, \cdot, \cdot\})$ be an LY-algebra. A linear map $D: T\rightarrow T$ is called a derivation of $T$, if it satisfies:
$$D([x, y])=[D(x), y]+[x, D(y)],$$
$$D([x, y, z])=[D(x), y, z]+[x, D(y), z]+[x, y, D(z)]$$
for all  $x, y, z\in T,$ that is to say, $D$ is simultaneously the derivation respect to the binary and ternary operation of $T$. The set of all derivations of $T$ is denoted by ${\rm Der}(T).$ ${\rm Der}(T)$  is a subalgebra of $\mathfrak{gl}(T)$ respect to the commutator operation.
\end{defi}

For an LY-algebra $(T, [\cdot, \cdot], \{\cdot, \cdot, \cdot\})$ and $x, y\in T,$ the linear map $\begin{array}[t]{rrl}
                                                                        L(x, y): T & \rightarrow & T \\
                                                                        z & \mapsto & \{x, y, z\}
                                                                      \end{array}
$ is (according to $(LY5)$ and $(LY6)$) a derivation of $T$. Note $L(T, T)=\{\sum L(x, y)| x, y\in T\}.$ By $(LY6)$, $L(T, T)$ is closed under commutation and $L(T, T)$ is a subalgebra of ${\rm Der}(T)$.

\begin{defi}[Generalized derivation, Quasi-derivation]
$D\in {\rm End}(T)$ is called a generalized derivation of $T$, if there exist $D^{(1)}, D^{(2)}, D^{(3)}, D^{(4)}, D^{(5)}\in {\rm End}(T)$ such that $$[D(x), y]+[x, D^{(1)}(y)]=D^{(2)}([x, y]),$$
$$\{D(x), y, z\}+\{x, D^{(3)}(y), z\}+\{x, y, D^{(4)}(z)\}=D^{(5)}([x, y, z])$$
for all $x, y, z\in T.$ $D$ is called a quai-derivation of $T$, if there exist $D', D''\in {\rm End}(T)$ such that $$[D(x), y]+[x, D(y)]=D'([x, y]),$$
$$\{D(x), y, z\}+\{x, D(y), z\}+\{x, y, D(z)\}=D''(\{x, y, z\})$$ for all $x, y, z\in T.$
We denote respectively by ${\rm GDer}(T)$ and ${\rm QDer}(T)$ the set of all generalized derivations and quasi-derivations of $T$. Both of them are subalgebras of $\mathfrak {gl}(T)$, and ${\rm QDer}(T)$ is a subalgebra of ${\rm GDer}(T)$.
\end{defi}
\begin{remark}
It is easy to see that a generalized derivation of an LY-algebra preserve its center.
\end{remark}
\begin{defi}[Centroid]
The set $$C(T)=\{D\in {\rm End}(T)|[D(x), y]=D([x, y]), \, {\mbox and }\, \{D(x), y, z\}=D(\{x, y, z\}), \forall x, y, z\in T\}$$ is called the centroid of $T.$
\end{defi}
$C(T)$ is closed under composition.  It is easy to show that, for a centerless LY-algebra $T$, $C(T)$ is commutative.
\begin{remark}
$C(T)$ is a subalgebra of $\mathfrak{gl}(T)$.
\end{remark}
\begin{defi}[Quasi-centroid]
The set $QC(T)=\{D\in {\rm End}(T)|[D(x), y]=[x, D(y)],$  and $\{D(x), y, z\}=\{x, D(y), z\}=\{x, y, D(z)\}, \forall x, y, z\in T\}$ is called the quasi-centroid of $T.$
\end{defi}
\begin{remark}\label{rem1}
\begin{enumerate}
\item It is easy to verify that $C(T)\subseteq {\rm QDer}(T)\cap QC(T).$
\item Although ${\rm Der}(T)$ and $C(T)$ preserve the derived algebra of $T$, neither ${\rm QDer}(T)$ nor $QC(T)$ need do so.
\end{enumerate}
\end{remark}
\begin{exam}
 Let $T$ be a two-dimensional LY-algebra spanned by $x, y$ with $[x, y]=y, \{x, y, y\}=y, \{y, x, x\}=0.$ The linear map $D: x\mapsto 0, y\mapsto x$ is a quasi-derivation, but $D([T, T]+\{T, T,T\})\nsubseteq [T, T]+\{T, T,T\}.$
\end{exam}
\begin{exam}
Let $T$ have a basis $x_0, x_1, \cdots, x_5$ with $$[x_0, x_1]=x_1, [x_0, x_3]=x_3, [x_0, x_5]=x_5, [x_1, x_2]=x_5, [x_3, x_4]=x_5,$$ $$\{x_0, x_1, x_0\}=x_1, \{x_0, x_3, x_0\}=x_3, \{x_0, x_1, x_1\}=\{x_0, x_1, x_3\}=\{x_0, x_3, x_1\}=\{x_3, x_1, x_1\}=x_5, $$ $$\{x_1, x_3, x_3\}=\{x_0, x_3, x_3\}=\{x_1, x_2, x_0\}=\{x_0, x_1, x_2\}=\{x_3, x_4, x_0\}=\{x_3, x_0, x_4\}=x_5,$$ and with other products $0$. Then $C(T)$ is spanned by ${\rm id}_T$ and $f_1, f_2$, where $f_1(x_0)=x_2, f_1(x_1)=-x_5, f_2(x_0)=x_4, f_2(x_3)=-x_5,$ while otherwise, $f_i(x_j)=0.$ And $QC(T)$ is spanned by $C(T)$ and $f_3$, where $f_3(x_1)=-x_4, f_3(x_3)=x_2$ with $f_3(x_i)=0$ for $i\neq 1, 3.$ Then $T^{(1)}$ is spanned by $x_1, x_3, x_5$, but $f_3(T^{(1)})\nsubseteq T^{(1)}$.
\end{exam}
\begin{defi}[central derivation]
A linear map $D: T\rightarrow T$ is called a central derivation of $T$ if it satisfies $[D(x), y]=D([x, y])=0, \, {\rm and}\, \{D(x), y, z\}=D(\{x, y, z\})=0, \forall x, y, z\in T. $\\
 The set of all central derivations of $T$ is denoted by ${\rm ZDer}(T).$ ${\rm ZDer}(T)$  is an ideal of $\mathfrak {gl}(T)$. It is clear that ${\rm ZDer}(T)\subseteq C(T).$
\end{defi}
It is easy to verify that ${\rm ZDer}(T)\subseteq {\rm Der}(T)\subseteq {\rm QDer}(T)\subseteq {\rm GDer}(T)\subseteq \mathfrak{gl}(T).$

For convenience, we use the following notations in section 5 and Lemma \ref{lem3} :\\
Suppose $T$ is a finite-dimensional LY-algebra with multiplications $\mu_1: T\times T\rightarrow T$ and $\mu_2: T\times T\times T\rightarrow T.$ Let $\Delta(T)$ denote the set of $(f, f^{(1)}, f^{(2)}, f^{(3)}, f^{(4)}, f^{(5)})\in {\rm End}(T)^6$ such that $\forall x, y, z\in T,$ $$\mu_1((f(x), y)+(x, f^{(1)}(y)))=f^{(2)}\circ \mu_1(x, y),$$ and
 \quad $\mu_2((f(x),y, z)+(x, f^{(3)}(y), z)+(x, y, f^{(4)}(z)))=f^{(5)}\circ \mu_2(x, y, z).$
 Then \Bea{\rm GDer}(T)&=&\{f\in {\rm End}(T)|\exists f^{(1)}, \cdots, f^{(5)}: (f, f^{(1)}, \cdots, f^{(5)})\in \Delta(T)\},\\
 {\rm Der}(T)&=& \{f\in {\rm End}(T)| (f, f, f, f, f, f)\in \Delta(T) \},\\
 C(T)&=& \{f\in {\rm End}(T)| (f, 0, f, 0, 0, f)\in \Delta(T)\},\\
 {\rm QDer}(T)&=& \{f\in {\rm End}(T)|\exists f', f''\in {\rm End}(T): (f, f, f', f, f, f'')\in \Delta(T)\},\\
 QC(T)&=&\{f\in {\rm End}(T)|(f, -f, 0, -f, 0, 0)\in \Delta(T)\, \mbox {and}\, (f, -f, 0, 0, -f, 0)\in \Delta(T)\}.\Eea
 Notice that if $(f, f^{(1)}, f^{(2)}, f^{(1)}, f^{(4)}, f^{(5)})\in \Delta(T)$, then $(f^{(1)}, f, f^{(2)}, f, f^{(4)}, f^{(5)})\in \Delta(T).$

\section{General results}
\setcounter{equation}{0}
\begin{lem}\label{lem1}
Let $T$ be an LY-algebra, then
\begin{enumerate}
\item [(1)] $[{\rm Der}(T), C(T)]\subseteq C(T).$
\item [(2)]$[{\rm QDer}(T), QC(T)]\subseteq QC(T).$
\item [(3)]$C(T)\cdot {\rm Der}(T)\subseteq {\rm Der}(T).$
\item [(4)] $C(T)\subseteq {\rm QDer}(T).$
\item [(5)] $[QC(T), QC(T)]\subseteq {\rm QDer}(T).$
\item [(6)]${\rm QDer}(T)+QC(T)\subseteq {\rm GDer}(T).$
\end{enumerate}
\end{lem}
\proof
(1)-(5) are easy. \\
(6) For ${\rm QDer}(T)+QC(T) \subseteq {\rm GDer}(T):$\\
Let $D_1\in {\rm QDer}(T), D_2\in QC(T).$ Then $\exists D'_1, D''_1\in {\rm End}(T), \forall x, y, z\in T,$ $$[D_1(x), y]+[x, D_1(y)]=D'_1([x, y]),$$
$$\{D_1(x), y, z\}+\{x, D_1(y), z\}+\{x, y, D_1(z)\}=D''_1(\{x, y, z\}).$$
Thus $\forall x, y, z\in T,$
\Bea
[(D_1+D_2)(x),y]&=& [D_1(x), y]+[D_2(x), y]\\
&=&D'_1([x, y])-[x, D_1(y)]+[x, D_2(y)]\\
&=&D'_1([x, y])-[x, (D_1-D_2)(y)]
\Eea
\Bea
\{(D_1+D_2)(x), y, z\}&=&\{D_1(x), y, z\}+\{D_2(x), y, z\}\\
&=&D''_1(\{x, y, z\})-\{x, D_1(y), z\}-\{x, y, D_1(z)\}+\{x, D_2(y), z\}\\
&=&D''_1(\{x, y, z\})-\{x, (D_1-D_2)(y), z\}-\{x, y, D_1(z)\}.
\Eea
Therefore, $D_1+D_2\in {\rm GDer}(T).$
\endproof

\begin{prop}
If $T$ is an LY-algebra, then $QC(T)+[QC(T), QC(T)]$ is a subalgebra of ${\rm GDer}(T).$
\end{prop}
\proof
By Lemma\ref{lem1}, (5), (6), we have $$QC(T)+[QC(T), QC(T)]\subseteq {\rm GDer}(T),$$
and \Bea
&&[QC(T)+[QC(T), QC(T)], QC(T)+[QC(T),QC(T)] ]\\
&\subseteq&[QC(T)+{\rm QDer}(T), QC(T)+[QC(T),QC(T)]]\\
&\subseteq&[QC(T), QC(T)]+[QC(T), [QC(T), QC(T)]]+[{\rm QDer}(T), QC(T)]\\
&&+[{\rm QDer}(T), [QC(T), QC(T)]].
\Eea
It is easy to verify that $[{\rm QDer}(T), [QC(T), QC(T)]]\subseteq [QC(T), QC(T)]$ by the Jacobi identity of Lie algebra. Thus $$[QC(T)+[QC(T), QC(T)], QC(T)+[QC(T), QC(T)]]\subseteq QC(T)+[QC(T), QC(T)].$$
\endproof
\begin{prop}
If $T$ is an LY-algebra, then $[C(T), QC(T)]\subseteq {\rm Hom}(T, Z(T)).$ Moreover, if $Z(T)=\{0\}, $ then $[C(T), QC(T)]=\{0\}$.
\end{prop}
\proof
Let $D_1\in C(T), D_2\in QC(T),$ then for all $x, y, z\in T,$ we have
\Bea
[[D_1, D_2](x), y]&=&[D_1D_2(x), y]-[D_2D_1(x), y]\\
&=&D_1([D_2(x), y])-[D_1(x), D_2(y)]\\
&=&D_1([D_2(x), y]-[x, D_2(y)])\\
&=&0,
\Eea
and \Bea
\{[D_1, D_2](x), y, z\}&=&\{D_1D_2(x), y, z\}-\{D_2D_1(x), y, z\}\\
&=&D_1(\{D_2(x), y, z\})-\{D_1(x), D_2(y), z\}\\
&=&D_1(\{D_2(x), y, z\}-\{x, D_2(y), z\})\\
&=&0.
\Eea
\endproof

The following lemma gives a condition for getting an equation more accurate than the inclusion in Remark \ref{rem1} for centerless LY-algebras:
\begin{lem}\label{lem3}
Let  ${\rm char}\mathbb K=0, D\in {\rm End}(T)$. If $(D, D, D', D, D, \frac{3D'}2),(D, -D, 0, -D, 0, 0),$ $ (D, -D, 0, 0, -D, 0)\in \Delta(T),$ then $\forall x, y, z, u, v\in T,$
$$[x, D([y, z])]=[x, [D(y), z]],$$
$$\{D(\{u, v, y\}), x, z\}=\{\{D(u), v, y\}, x, z\}.$$
\end{lem}
\proof
Suppose $(D, D, D', D, D, \frac{3D'}2),(D, -D, 0, -D, 0, 0),$ $ (D, -D, 0, 0, -D, 0)\in \Delta(T),$ then for all $x, y, z\in T,$ $$[D(x), y]+[x, D(y)]=D'([x, y]),$$
$$ \{D(x), y, z\}+\{x, D(y), z\}+\{x, y, D(z)\}=D''(\{x, y, z\}),$$ and $$[D(x), y]=[x, D(y)],$$
$$\{D(x), y, z\}=\{x, D(y), z\}=\{x, y, D(z)\}.$$
Put $g=\frac{D'}2,$ then we have $\forall x, y, z\in T,$
\Bea
g(\{x, y, z\})&=&\frac 13(\{D(x), y, z\}+\{x, D(y), z\}+\{x, y, D(z)\})\\
&=&\{D(x), y, z\}=\{x, D(y), z\}=\{x, y, D(z)\}.
\Eea
We notice that $\forall u, v, x, y, z\in T,$
\Bea
&&\{u, v, g(\{x, y, z\})\}=\{u, v, \{D(x), y, z\}\}\\
&=& \{\{u, v, D(x)\}, y, z\}+\{D(x), \{u, v, y\}, z\}+\{D(x), y, \{u, v, z\}\}\\
&=&\{g(\{u, v, x\}), y, z\}+\{x, D(\{u, v, y\}), z\}+g(\{x, y, \{u, v, z\}\}).
\Eea
So, $\forall u, v, x, y, z\in T,$
\bea
&&\{D(\{u, v, y\}), x, z\}=g(\{\{u, v, y\}, x, z\})\nonumber\\
&=&\{g(\{u, v, x\}), y, z\}+g(\{x, y, \{u, v, z\}\})-\{u, v, g(\{x, y, z\})\},
\eea
\bea
&&g(\{y, \{u, v, x\}, z\})=-g(\{\{u, v, x\}, y, z\})\nonumber\\
&=&-\{g(\{u, v, y\}), x, z\}+\{u, v, g(\{y, x, z\})\}-g(\{y, x, \{u, v, z\}\}),
\eea

\bea
&&g(\{u, v, \{y, x, z\}\})\nonumber\\
&\overset{{\rm exchange} (u, v) {\rm and} (y,x)}{\underset{ {\rm in}(3.2)}=}&-g(\{u, \{y,x,v\}, z\})-g(\{\{y, x, u\}, v, z\})+\{y, x, g(\{u, v, z\})\}
\eea
Thus, we have $\forall u, v, x, y, z\in T,$
\bea
&&g(\{\{u, v, y\}, x, z\})=\{D(\{u, v, y\}), x, z\}\nonumber\\
&=&g(\{u, v, \{y, x, z\}\})-g(\{y, \{u, v, x\}, z\})-g(\{y, x, \{u, v, z\}\})\nonumber\\
&=&-g(\{u, \{y, x, v\}, z\})-g(\{\{y, x, u\}, v, z\})+\{y, x, g(\{u, v, z\})\}\nonumber\\
&&+\{g(\{u, v, y\}), x, z\}-\{u, v, g(\{y, x, z\})\}+g(\{y, x, \{u, v, z\}\})-g(\{y, x, \{u, v, z\}\})\nonumber\\
&=& -\{u, \{y, x, v\}, D(z)\}-\{\{y, x, u\}, v, D(z)\}+\{y, x, \{u, v, D(z)\}\}\nonumber\\
&&+\{g(\{u, v, y\}), x, z\}-\{u, v, \{y, x, D(z)\}\}\nonumber\\
&=&\{g(\{u, v, y\}), x, z\}=\{\{D(u), v, y\}, x, z\}.
\eea
On the other hand, we have \begin{equation}\label{eq1}
g([x, y])=\frac 12([D(x), y]+[x, D(y)])=[D(x), y]=[x, D(y)]
\end{equation}
We notice that $\forall x, y, z\in T,$
\Bea
[g([x, y]), z]&=& [[D(x), y], z]\\
&=& [[D(x), z], y]+[D(x), [y, z]]+\{D(x), z, y\}+\{z, y, D(x)\}+\{y, D(x), z\}\\
&=&[g([x, z]), y]+[x, D([x, z])]+g(\{x, z, y\}+\{z, y, x\}+\{y, x, z\})
\Eea
So that, $\forall x, y, z\in T,$
\bea
[g([x, y]), z]+[g([z, x]), y]&=&[x, D([y, z])]+g(\{x, z, y\}+\{z, y, x\}+\{y, x, z\})\label{eq5}\\
&=&g([x, [y, z]])+g(\{x, z, y\}+\{z, y, x\}+\{y, x, z\})\label{eq2}
\eea
Adding the three equations obtained from Eq.(\ref{eq2}) by cyclically permuting $x, y, z$, %and using Eq. (\ref{eq2}),
we get $\forall x, y, z\in T,$
\Bea &&2([g([x, y]), z]+[g([z, x]), y]+[g([y, z]), x])\\
&=&g([x, [y, z]]+[y, [z, x]]+[z, [x, y]])+3g(\{x, z, y\}+\{z, y, x\}+\{y, x, z\})\\
&=&g([x, [y, z]]+[y, [z, x]]+[z, [x, y]]+\{x, z, y\}+\{z, y, x\}+\{y, x, z\})\\
&=&2g(\{x, z, y\}+\{z, y, x\}+\{y, x, z\}).
\Eea
Thus, $\forall x, y, z\in T,$ \bea
&&[g([x, y]), z]+[g([z, x]), y]\nonumber\\
&=&-[g([y, z]), x]+g(\{x, z, y\}+\{z, y, x\}+\{y, x, z\})\nonumber\\
&=&[x, g([y, z])]+g(\{x, z, y\}+\{z, y, x\}+\{y, x, z\})\label{eq3}.\eea
Comparing Eq. (\ref{eq3}) with Eq.(\ref{eq5}), we have $\forall x, y, z\in T, [x, D([y, z])]=[x, g([y, z])].$ By Eq. (\ref{eq1}), $\forall x, y, z\in T, $
\begin{equation}\label{eq4}
[x, D([y, z])]=[x, [D(y), z]]
\end{equation}
%Since $Z(T)=\{0\},$ from (3.7) and (3.8), we have $\forall u, v, y\in T,$ $$D(\{u, v, y\})=g(\{u, v, y\})=\{D(u), v, y\}=\{u, D(v), y\}=\{u, v, D(y)\}$$
%and\qquad\qquad\qquad  $\forall x, y, z\in T, D([y, z])=f([y, z])=[D(y), z]=[y, D(z)].$ \\
%Therefore,  $D\in C(T).$
\endproof\\
Then we have the following.
\begin{prop}
Let ${\rm char}\mathbb K=0, S=\{D\in {\rm End}(T)|\exists D'\in {\rm End}(T):(D, D, D', D, D, \frac{3D'}2)\in \Delta(T)\}$ be a subset of ${\rm QDer}(T)$. If $Z(T)=\{0\},$ then $C(T)=S\cap QC(T).$
\end{prop}
Recall the definition of Jordan algebra:
\begin{defi}
Let $L$ be an algebra over $\mathbb K$. If the multiplication satisfies the following identities: $$x y= y x, \quad (xy)(xx)=x(y(xx)),$$ for all $x, y\in L$, then we call $L$ a Jordan algebra.
\end{defi}
\begin{prop}
Let $T$ be an LY-algebra over $\mathbb K$, then ${\rm End}(T)$ is a Jordan algebra with the multiplication $\begin{array}[t]{rrl}
                  \ast:  {\rm End}(T) \times {\rm End}(T)& \rightarrow  & {\rm End}(T) \\
                  (D_1, D_2) & \mapsto & D_1D_2+D_2D_1
                \end{array}.
$
\end{prop}
\begin{coro}\label{coro1}
Let $T$ be an LY-algebra over $\mathbb K$, then $QC(T)$ is a Jordan algebra with the multiplication $\begin{array}[t]{rrl}
                  \star:  QC(T) \times QC(T)& \rightarrow  & QC(T) \\
                  (D_1, D_2) & \mapsto & D_1D_2+D_2D_1
                \end{array}.$
\end{coro}
\begin{thm}
Let  $T$ be an LY-algebra over $\mathbb K$. We have
\begin{enumerate}
\item [(1)] If ${\rm char}\mathbb K\neq 2,$ then $QC(T)$ is a Lie algebra with commutator if $QC(T)$ is also an associative algebra with respect to composition.
\item [(2)] If ${\rm char}\mathbb K\notin \{2, 3\},$ and $Z(T)=\{0\}$, then $QC(T)$ is a Lie algebra if and only if $$[QC(T), QC(T)]=\{0\}.$$
\end{enumerate}
\end{thm}
\proof
\begin{enumerate}
\item [(1)] $``\Longleftarrow":$ For all $D_1, D_2\in QC(T),$ we have
$D_1D_2\in QC(T) $ and  $D_2D_1\in QC(T).$\\
So, $[D_1, D_2]=D_1D_2-D_2D_1\in QC(T).$
 Hence, $QC(T)$ is a Lie algebra.\\
$``\Longrightarrow":$ Let $D_1, D_2\in QC(T).$ Suppose $[D_1, D_2]\in QC(T)$. \\
Notice that $D_1D_2=D_1\star D_2+\frac {[D_1, D_2]}2.$ By Corollary \ref{coro1} we have $D_1\star D_2\in QC(T).$\\ So, $D_1D_2\in QC(T).$
\item [(2)] $``\Longrightarrow":$ Suppose that $D_1D_2\in QC(T).$\\
For all $x, y, z\in T, QC(T)$ is a Lie algebra, so $[D_1, D_2]\in QC(T).$ Then $[[D_1, D_2](x), y]=[x, [D_1, D_2](y)].$ On the other hand, by Lemma \ref{lem1}(5) we have \bea 2[[D_1, D_2](x), y]=0, \forall x, y\in T.\eea and \bea
3\{[D_1, D_2](x), y, z\}=0, \forall x, y, z\in T.
\eea
Since ${\rm Char}\mathbb K\notin \{2, 3\}, $ $[[D_1, D_2](x), y]=0$  and $\{[D_1, D_2](x), y, z\}=0, \forall x, y, z\in T.$ What's more, $Z(T)=\{0\},$ so, $[D_1, D_2]=0.$\\
$``\Longleftarrow":$ Trivial.
\end{enumerate}
\endproof

\begin{lem}$^{[\ref{ref19}]}$
Let $V$ be a vector space and $f: V\rightarrow V$ a linear map. $\pi_f$ denotes the minimal polynomial of $f$. If $X^2\nmid \pi_f, $ then $V={\rm Ker}(f)\overset{\cdot}+{\rm Im}(f)$.
\end{lem}
Similar to [\ref{ref19}], we can prove the following results:
\begin{prop}
Let $T$ be an LY-algebra, $D\in C(T)$. Then
\begin{enumerate}
\item [(1)] ${\rm Ker}(D)$ and ${\rm Im}(D)$ are ideals of $T$.
\item [(2)] Suppose $T$ is indecomposable, $D\in C(T)\setminus\{0\}$. If $X^2|\pi_D$, then $D$ is invertible.
\item [(3)] Suppose $T$ is perfect. If $T$ is indecomposable and $C(T)$ consists of semisimple elements, then $C(T)$ is a field.
\end{enumerate}
\end{prop}

\begin{lem}
Let $T$ be a centerless LY-algebra. If $D\in QC(T)$ and $X^3|\nmid \pi_D,$ then $T={\rm Ker}(D)\oplus {\rm Im}(D).$
\end{lem}

\begin{coro}
Let $T$ be a centerless and indecomposable LY-algebra over an algebraically closed field $\mathbb K$. If $D\in QC(T)$ is semisimple, then $D\in Z_{C(T)}(GDer(T))$.
\end{coro}

\section{Quasi-derivation embedded as derivation of a larger LY-algebra}
Inspired by [\ref{ref20}], the quasi-derivations of an LY-algebra can be embedded as derivations in a larger LY-algebra.  For this, let $T$ be an LY-algebra over $\mathbb K$ and $t$ an indeterminant. We define
$$ \check{T}=\bigoplus\limits_{1\leq i\leq 3}(T\otimes t^i),$$ and multiplications on it:
$$[xt^i, yt^j]=\left\{
                 \begin{array}{ll}
                   [x, y]t^{i+j}, & \mbox{if}\, i+j=2 \\
                   0, & \mbox{if}\, i+j>2
                 \end{array}
               \right.
$$
$$\{xt^i, yt^j, zt^k\}=\left\{
                 \begin{array}{ll}
                   \{x, y, z\}t^{i+j+k}, & \mbox{if}\, i+j+k=3 \\
                   0, & \mbox{if}\, i+j+k>3
                 \end{array}
               \right.$$
Then $(\check{T}, [\cdot, \cdot], \{\cdot, \cdot, \cdot\})$ is an LY-algebra.

If $U, V$ are two subspace of $T$ satisfies $T=U\overset{\cdot}+[T, T]=V\overset{\cdot}+\{T, T, T\},$ then $$\check{T}=Tt+Tt^2+Tt^3=Tt+Ut^2+[T, T]t^2+Vt^3+\{T, T, T\}t^3.$$
Define $\varphi: QDer(T)\rightarrow {\rm End}(\check{T})$ as follow:\\
for $(D, D, D', D, D, D'')\in \Delta (T),$
$$\varphi(D)(at+ut^2+bt^2+vt^3+ct^3)=D(a)t+D'(b)t^2+D''(c)t^3,$$
$\forall a\in T, u, v \in U, b\in [T, T], c\in\{T, T, T\}.$
\begin{prop}\label{prop1}
$T, \check{T}, \varphi$ are defined as above. Then
\begin{enumerate}
\item [(1)] $\varphi$ is injective and $\varphi(D)$ does not depend on the choice of the pair $(D', D'')$.
\item [(2)] $\varphi(QDer(T))\subseteq {\rm Der}(\check{T}).$
\end{enumerate}
\end{prop}
\proof
\begin{enumerate}
\item [(1)] Let $D_1, D_2\in QDer(T)$ such that $\varphi(D_1)\subseteq \varphi(D_2).$ Then for all $a\in T_1, b\in [T, T], c\in [T, T, T], u\in U, v\in V$, we have $$\varphi(D_1)(at+ut^2+bt^2+vt^3+ct^3)=\varphi(D_2)(at+ut^2+bt^2+vt^3+ct^3),$$
i.e. $D_1(a)t+D'_1(b)t^2+D''_1(c)t^3=D_2(a)t+D'_2(b)t^2+D''_2(c)t^3.$\\
So, we proved that $D_1(a)=D_2(a), \forall a\in T.$\\
Hence, $D_1=D_2$ and $\varphi$ is injective. Suppose that there exists $\widetilde{D'}$ and $\widetilde{D''}$ such that $(D, D, \widetilde{D'}, D, D, \widetilde{D''})$ is also in $\Delta(T).$ Then \Bea
\varphi(D)(at+ut^2+bt^2+vt^3+ct^3)&=&D(a)t+D'(b)t^2+D''(c)t^3\\
&=&D(a)t+\widetilde{D'}(b)t^2+\widetilde{D''}(c)t^3,\Eea
and\qquad \, $D'([x, y])=[D(x), y]+[x, D(y)]=\widetilde{D'}([x, y]),$ $$D''\{x, y, z\}=\{D(x), y, z\}+\{x, D(y),z\}+\{x, y, D(z)\}=\widetilde{D''}(\{x, y, z\})$$
for all $x, y, z\in T.$ Thus, $D'(b)=\widetilde{D'}(b)$ and $D''(c)=\widetilde{D''}(c).$ Therefore, $\varphi(D)$ does not depend on the choice of $(D', D'')$.
\item[(2)]Let $D\in QDer(T)$. By definition, for all $i+j\geq 3, [xt^i, yt^j]=0$ and for all $i+j+k\geq 4, \{xt^i, yt^j, zt^k\}=0.$ Thus, to show $\varphi(D)\in {\rm Der}(\check{T}),$ it suffice to check the following equations: $$\varphi(D)([xt, yt])=[\varphi(D)(xt), yt]+[xt, \varphi(D)(yt)]$$
and $\varphi(D)\{xt, yt, zt\}=\{\varphi(D)(xt), yt, zt\}+\{xt, \varphi(D)(yt), zt\}+\{xt, yt, \varphi(D)(zt)\}$ for all $x, y, z\in T.$ In fact, for all $x, y, z\in T,$ we have
\Bea
\varphi(D)([xt, yt])&=&\varphi(D)([x, y]t^2)
=D'([x, y])t^2\\
&=&([D(x), y]+[x, D(y)])t^2\\
&=&[D(x)t, yt]+[xt, D(y)t]\\
&=&[\varphi(D)(xt), yt]+[xt, \varphi(D)(yt)],
\Eea
\Bea
\varphi(D)(\{xt, yt, zt\})&=&\varphi(D)\{x, y, z\}t^3=D''(\{x, y, z\})t^3\\
&=&(\{D(x), y, z\}+\{x, D(y), z\}+\{x, y, D(z)\})t^3\\
&=&\{D(x)t, yt, zt\}+\{xt, D(y)t, zt\}+\{xt, yt, D(z)t\}\\
&=&\{\varphi(D)(xt), yt, zt\}+\{xt, \varphi(D)(yt), zt\}+\{xt, yt, \varphi(D)(zt)\}.
\Eea
Therefore,  $\varphi(QDer(T))\subseteq {\rm Der}(T).$
\end{enumerate}
\endproof
\begin{prop}\label{prop2}
If $T$ is a centerless LY-algebra and $\check{T}, \varphi$ are defined as above, then $${\rm Der}(\check{T})=\varphi(QDer(T))\oplus ZDer(\check{T}).$$
\end{prop}
\proof
Suppose that $xt+yt^2+zt^3\in Z(\check{T})$, then $\forall u_i, v_i, w_i\in T, i\in\{1, 2\}$, $$0=[xt+yt^2+zt^3, u_1t+v_1t^2+w_1t^3]=[x, u_1]t^2$$ and \qquad\qquad \qquad $0=\{xt+yt^2+zt^3, u_1t+v_1t^2+w_1t^3\}=\{x, u_1, u_2\}t^3.$\\
 So that $[x, u_1]=0$ and $\{x, u_1, u_2\}=0,$ then $x\in Z(T).$ Since $Z(T)=\{0\},$ we have $x=0.$\linebreak Thus $Z(\check{T})\subseteq Tt^2\overset{\cdot}+Tt^3.$ The anti-conclusion is trivial.

We know that a derivation of a LY-algebra preserve its center: $$\forall g\in {\rm Der}(\check{T}), g(Z(\check{T}))\subseteq Z(\check{T}).$$
Hence $g(Ut^2+Ut^3)\subseteq g(Z(\check{T}))\subseteq Z(\check{T})=Tt^2\overset{\cdot}+Tt^3.$

Any linear map $f: Tt+Ut^2+Vt^3\rightarrow Tt^2+Tt^3$ extends to an element of $ZDer(\check{T})$ by taking $f([T, T]t^2+\{T, T, T\}t^3)=0.$ Thus, given any $g\in {\rm Der}(\check{T})$, we can define:
 $$\begin{array}{rll}
    f: &Tt+Ut^2+[T, T]t^2+Vt^3+\{T, T, T\}t^3  \rightarrow   Tt^2+Tt^3 \\
   & x  \mapsto   \left\{
                      \begin{array}{ll}
                        0, & \mbox{if}\, x\in [T, T]t^2+\{T, T, T\}t^3 \\
                        f(x) \, \mbox{such}\, \mbox{that}\, (g-f)(x)\in Tt, & \mbox{if}\, x\in Tt \\
                        g(x), & \mbox{if}\, x\in Ut^2+Vt^3
                      \end{array}
                    \right.
  \end{array}
$$
Then $f\in ZDer(\check{T})$ and $(g-f)(Tt)\subseteq Tt, (g-f)(Ut^2+Vt^3)=0.$ In addition, since\linebreak $\check{T}^{(1)}=[T, T]t^2+\{T, T, T\}t^3, $ $(g-f)([T, T]t^2+\{T, T, T\}t^3)\subseteq [T, T]t^2+\{T, T, T\}t^3.$ Thus, there exist $D, D', D''\in {\rm End}(T)$ such that $\forall a\in T, \forall b\in [T, T], \forall c\in \{T, T, T\},$ $$(g-f)(at)=D(a)t, (g-f)(bt^2)=D'(b)t^2, (g-f)(ct^3)=D''(c)t^3.$$ Since $g-f\in {\rm Der}(\check{T})$ (for $ZDer(T)\subseteq {\rm Der}(\check{T})$) and by the definition of ${\rm Der}(\check{T})$, we have $\forall a_1, a_2, a_3, a_4, a_5\in T,$ $$[(g-f)(a_1t), a_2t]+[a_1t, (g-f)(a_2t)]=(g-f)[a_1t, a_2t],$$
and $\{(g-f)(a_3)t, a_4t, a_5t\}+\{a_3t, (g-f)(a_4t), a_5t\}+\{a_3t, a_4t,(g-f)(a_5t)\}=(g-f)\{a_3t, a_4t, a_5t\}.$ Hence, $$[D(a_1),a_2]+[a_1, D(a_2)]=D'([a_1, a_2]),$$
$$\{D(a_3), a_4, a_5\}+\{a_3, D(a_4), a_5\}+\{a_3, a_4, D(a_5)\}=D''(\{a_3, a_4, a_5\}).$$ Thus, $(D, D, D', D, D, D'')\in \Delta(T).$ Therefore, $g-f=\varphi(D)\in \varphi(QDer(T))$. So, ${\rm Der}(\check{T})\subseteq \varphi(QDer(T))+ZDer(T).$ From Proposition \ref{prop1}(2), we konw ${\rm Der}(\check {T})=\varphi(QDer(T))+ZDer(T).$ Now, we prove that $\varphi(QDer(T))\cap ZDer(T)=\{0\}.$ In fact, $\forall f\in \varphi(QDer(T))\cap ZDer(T), \exists D\in QDer(T): f=\varphi(D)$ and $f\in ZDer(\check{T}).$ Then
\Bea
f(at+ut^2+bt^2+vt^3+ct^3)
&=& \varphi(D)(at+ut^2+bt^2+vt^3+ct^3)\\
&=&D(a)t+D'(b)t^2+D''(c)t^3,
\Eea and \qquad \qquad $f(at+ut^2+bt^2+vt^3+ct^3)\in Z(\check{T})=Tt^2+Tt^3,$\\
for all $a\in T, b\in [T, T], c\in \{T, T, T\}, u\in U, v\in V.$  That is to say, $D(a)=0, \forall a\in T.$ So, $D=0.$ It follows that $f=0.$ Thus, ${\rm Der}(\check{T})= \varphi(QDer(T))\overset{\cdot}+ZDer(T).$
\endproof

Observe that, in the case of Proposition \ref{prop2}, $\varphi(QDer(T))$ may be viewed as the natural image of ${\rm Der}(\check{T})$ in ${\rm End}(\check{T}/\check{T}^{(1)}).$
%\section{LY-algebras with ${\rm End}(T)=QDer(T)$}
%For an LY-algebra $T,$ we may consider
\section{Quasi-derivations and Robustness}
In this section, we suppose $\mathbb K$ is a field with characteristic zero.
\begin{defi} Let $(T, \mu_1, \mu_2)$ be an LY-algebra. If $f$ a nonsingular element of ${\rm End}(T)$ such that $(T, f\circ \mu_1, f\circ \mu_2)$ is an LY-algebra, we call $(T, f\circ \mu_1, f\circ \mu_2)$ a perturbation of $(T, \mu_1, \mu_2)$. The perturbation is said to be inessential if $f\circ \mu_i=c\circ \mu_i,i=1,2,$ for some $c\in C(T)$. We say $(T, \mu_1, \mu_2)$ is robust if every perturbation of $(T, \mu_1, \mu_2)$ is inessential.
\end{defi}

Let $(T, \mu_1, \mu_2)$ be an LY-algebra and $V$ be a $T-$module. $(\rho, D, \theta; V)$ denotes the corresponding representation of $T$, which is also abbreviated as $(\rho, \theta)$.  If $(\rho, \theta)=(0, 0)$, we call $V$ a trivial $T-$module. Here we use the regular representation. $Z^{2p}(T, V)\times Z^{2p+1}(T, V), B^{2p}(T, V)\times B^{2p+1}(T, V)$ and $H^{2p}(T, V)\times H^{2p+1}(T, V)$ denote the $p-$th cocycles, coboundaries and cohomologies of $T$ with coefficients in $V$, while $\delta=(\delta_{\rm I}, \delta_{\rm II})$ denotes the coboundary operator. For each $(f, g)\in C^{2}(T, V)\times C^{3}(T, V)$ there is  another coboundary operation $\delta^{*}=(\delta^{*}_{\textrm{I}}, \delta^{*}_{\textrm{II}})$ of $C^{2}(T, V)\times C^{3}(T, V)$ into $C^{3}(T, V)\times C^{4}(T, V)$. For more details, see [\ref{ref091}]. $\overset{\circ}T$ denotes the trivial $T-$module on the underlying vector space of $T$, while the regular module is still denoted by $T$. To distinguish the coboundary operator $\delta$ on $C^{2p}(T, T)\times C^{2p+1}(T, T)$, we denote $\overset{\circ}{\delta}$ the coboundary map on $C^{2p}(T, \overset{\circ}T)\times C^{2p+1}(T, \overset{\circ}T)$. If we define $B^2(T,\overset{\circ}T)\times B^3(T, \overset{\circ}T)=\{\overset{\circ}{\delta}(f, g)|f,g\in C^1(T, \overset{\circ}T)\}$, then we observe that $B^2(T, \overset{\circ}T)\times B^3(T, \overset{\circ}T)=({\rm End}(T)\circ \mu_1)\times({\rm End}(T)\circ \mu_2)$. In addition, if $f\in {\rm End}(T)$, then $\delta(f\circ \mu_1, f\circ \mu_2)=\overset{\circ}{\delta}(f\circ \mu_1, f\circ \mu_2).$

Followed by Yamaguti, $H^{1}(T, T)=\{f\in C^1(T, T)|\delta_{\rm I}(f)=0=\delta_{\rm II}(f)\}={\rm Der}(T).$

After a simple calculation, we get
\begin{prop}\label{prop4}
Let $(T, \mu_1, \mu_2)$ be an LY-algebra. If $f\in {\rm End}(T)$ then \begin{center}
$(T, f\circ \mu_1, f\circ\mu_2)$ is an LY-algebra $\Longleftrightarrow$ $(f\circ \mu_1, 2f\circ \mu_2)\in Z^2(T, T)\times Z^3(T, T).$\end{center}
\end{prop}

This leads to the following:
\begin{coro}
Let $(T, \mu_1, \mu_2)$ be an LY-algebra. If $(Z^2(T, T)\times Z^3(T, T))\cap(B^2(T, \overset{\circ}{T})\times B^3(T, \overset{\circ}{T}))=\{(c\circ\mu_1, 2c\circ\mu_2)|c\in C(T)\},$ then $T$ is robust.
\end{coro}

The following lemma reveals the relationship between quasi-derivation and cohomology.
\begin{lem}\label{lem2}
Let $(T, \mu_1, \mu_2)$ be an LY-algebra, $f, f', f''\in {\rm End}(T).$ Then
\begin{enumerate}
\item [(1)] $f\in {\rm QDer}(T)$ if and only if $\delta(f, f)\in B^2(T, \overset{\circ}T)\times B^3(T, \overset{\circ}T)$. More specifically,
$$(f, f, f',f,f,f'')\in \Delta(T) \Longleftrightarrow (\delta_{\rm I}(f), \delta_{\rm II}(f))=((f'-f)\circ \mu_1, (f''-f)\circ \mu_2).$$
\item [(2)] If $(f, f, f',f,f,f'')\in \Delta(T),$ then
\Bea
&&-\mu_1( x_3, (f''-f)\circ \mu_2(x_1, x_2, x_4))-\mu_1((f''-f)\circ\mu_2(x_1, x_2, x_3), x_4)\\
&&+(f''-f)\circ\mu_2(x_1, x_2, \mu_1(x_3, x_4))
+\mu_2(x_1, x_2, (f'-f)([x_3, x_4]))\\
&&-((f'-f)\circ \mu_1(\mu_2(x_1, x_2, x_3), x_4)-(f'-f)\circ\mu_1(x_3, \mu_2(x_1, x_2, x_4))=0
\Eea
and \Bea
&&-\mu_2((f''-f)\circ\mu_2(x_1, x_2, x_3), x_4, x_5)+\mu_2((f''-f)\circ \mu_2(x_1, x_2, x_4), x_3, x_5)\\
&&+\mu_2(x_1, x_2, (f''-f)\circ \mu_2(x_3, x_4, x_5)
-\mu_2(x_3, x_4, f\circ\mu_2(x_1, x_2, x_5))\\
&&-(f''-f)\circ\mu_2(\mu_2(x_1, x_2, x_3), x_4, x_5)-(f''-f)\circ\mu_2(x_3, \mu_2(x_1, x_2, x_4), x_5)\\
&&-(f''-f)\circ\mu_2(x_3, x_4,\mu_2(x_1, x_2, x_5))+(f''-f)\circ\mu_2(x_1, x_2, \mu_2(x_3, x_4, x_5))=0.
\Eea

\end{enumerate}
\end{lem}
\proof
(2) Suppose $(f, f, f',f,f,f'')\in \Delta(T).$ According to (1), $(\delta_{\rm I}(f), \delta_{\rm II}(f))=((f'-f)\circ \mu_1, (f''-f)\circ \mu_2)$. We konw that $\delta_{\textrm{I}}\delta_{\textrm{I}}f=\delta^{*}_{\textrm{I}}\delta_{\textrm{I}}f=0$ and $\delta_{\textrm{II}}\delta_{\textrm{II}}f=\delta^{*}_{\textrm{II}}\delta_{\textrm{II}}f=0$ from Yamaguti ([\ref{ref091}]). So  $(\delta_{\rm I}(f), \delta_{\rm II}(f)) \in Z^2(T, T)\times Z^3(T, T)$, where $$Z^2(T, T)\times Z^3(T, T)=\{(f, g)\in C^{2}(T, V)\times C^{3}(T, V)|\delta_{\textrm{I}}f=\delta^{*}_{\textrm{I}}f=0, \delta_{\textrm{II}}g=\delta^{*}_{\textrm{II}}g=0 \}.$$ By the definition of $\delta_{\rm I}$ and $\delta_{\rm II}$ we get the conclusion.
\endproof
%\begin{lem}
%$A=\{f\in {\rm End(T)}|\exists f'\in {\rm End}(T): (f, f, f', f, f, f, f')\in \Delta(T)\}$ is a subalgebra of the Lie algebra ${\rm QDer}(T)$.
%\end{lem}
\begin{prop}
 Let $(T, \mu_1, \mu_2)$ be an LY-algebra, then  ${\rm QDer}(T)={\rm Der}(T)+C(T)\Longrightarrow (B^2(T, T)\times B^3(T, T))\cap ( B^2(T, \overset{\circ}T)\times B^3(T, \overset{\circ}T))=\{(c\circ\mu_1, 2c\circ\mu_2)|c\in C(T)\}.$
\end{prop}
\proof Suppose ${\rm QDer}(T)={\rm Der}(T)+C(T).$ Let $f, g_1, g_2\in {\rm End}(T)$  such that $\delta(f, f)=\overset{\circ}{\delta}(g_1, g_2),$ then $(f, f, f+g_1, f, f, f, f+g_2)\in \Delta(T)$ by Lemma \ref{lem2}(1), and $f\in {\rm QDer}(T)$. In addition, we have $f=D+c$, with $D\in {\rm Der}(T), c\in C(T)$, so that $$\delta (f, f)=(\delta_{\rm I}(D), \delta_{\rm II}(D))+(\delta_{\rm I}(c), \delta_{\rm II}(c))=(\delta_{\rm I}(c), \delta_{\rm II}(c))=(c\circ \mu_1, 2c\circ \mu_2).$$
Thus, $(B^2(T, T)\times B^3(T, T))\cap ( B^2(T, \overset{\circ}T)\times B^3(T, \overset{\circ}T))\subseteq\{(c\circ\mu_1, 2c\circ\mu_2)|c\in C(T)\}.$\\
Conversely, Let $c\in C(T)$. Then according to the definition of $\delta_{\rm I}$ and $\delta_{\rm II},$  for $f\in {\rm End}(T)$, by putting $f'=f+c, f''=f+2c$, we have $(f, f, f', f, f,f, f'')\in \Delta(T)$ by Lemma \ref{lem2}(1). So, $f\in {\rm QDer}(T).$ By Lemma \ref{lem2} (1), $(c\circ \mu_1, 2c\circ \mu_2)=\delta(f, f)\in (B^2(T, T)\times B^3(T, T))\cap ( B^2(T, \overset{\circ}T)\times B^3(T, \overset{\circ}T)).$
%$\Longleftarrow):$ Suppose that $B^2(T, T)\times B^3(T, T)=(C(T)\circ \mu_1)\times (C(T)\circ \mu_2).$ \\
%If $f\in A,$ then by Proposition\ref{prop3}, $\delta(f, f)=((f'-f)\circ \mu_1, (f'-f)\circ \mu_2)\in B^2(T, T)\times B^3(T, T).$ So, $((f'-f)\circ \mu_1,  (f'-f)\circ \mu_2)\in (C(T)\circ \mu_1)\times (C(T)\circ \mu_2)$, and $f'-f\in C(T)$.
\endproof

\begin{coro}
Let $(T, \mu_1, \mu_2)$ be an LY-algebra. If $H^2(T, T)\times H^3(T, T)=\{0\}$ and ${\rm QDer}(T)={\rm Der}(T)+C(T)$, then $T$ is robust.
\end{coro}

From [\ref{ref16}], we know that if $H^2(T, T)\times H^3(T, T)=\{0\},$ then $T$ is rigid. The cohomological condition above indicates the existence of many algebras that are both robust and rigid. %However, the two properties are independent to each other.
\section{Acknowledgements}
 The first author gratefully acknowledges the support of  NSFC (No. 11501564) and Fundamental Research Funds for the Central Universities
(Grant No. 3122015L005).
 The second author gratefully acknowledges the support of  NSFC (No. 11801066) and Fundamental Research Funds for the Central Universities.
The third author gratefully acknowledges the support of  NNSF of China
(No. 11771069), NSF of Jilin province (No. 20170101048JC) and the project of jilin province department of education (No. JJKH20180005K). %The second author gratefully acknowledges the support of.

\end{document}